\documentclass[10pt]{amsart}
\usepackage[english]{babel}
\usepackage{amssymb}
\usepackage{amsmath}
\newcommand{\ignore}[1]{}

\newcommand{\J}{\mathcal{J}}
\newcommand{\E}{\mathcal{E}}

\newcommand{\A}{\mathcal{A}}

\newcommand{\N}{\mathcal{N}}
\newcommand{\Pe}{\mathcal{P}}

\newcommand{\sts}{\mathcal{O}_X}
\newcommand{\Pic}{{\rm Pic}}

\newcommand{\lb}{\mathcal{L}}

\newtheorem{lemma}{Lemma}
\newtheorem{theorem}{Theorem}
\newtheorem{proposition}{Proposition}
\newtheorem{corollary}{Corollary}

\theoremstyle{definition}

\newtheorem{example}{Example}

\newtheorem{remark}{Remark}



\title{Forms of higher degree permitting composition}

\author{S. Pumpl\"un}
\email{susanne.pumpluen@nottingham.ac.uk}
\address{School of Mathematics\\
University of Nottingham\\
University Park\\
Nottingham NG7 2RD\\
United Kingdom
 }

\begin{document}

\begin{abstract} Nondegenerate forms $N$ of degree $d$ on a unital nonassociative algebra $A$ over a ring $R$
which permit composition, i.e., satisfy $N(1)=1$ and
$N(xy)=N(x)N(y)$ for all $x,y$ in $A$, are studied. These forms were first classified by Schafer over fields of
 characteristic $0$ or $>d$.
We investigate cubic and quartic nondegenerate forms which permit composition over certain rings and curves.
 Classes of highly degenerate cubic forms $N$ over fields which permit composition are constructed.
\end{abstract}

\maketitle

\section{Introduction}

 Finite-dimensional forms of degree $d$ permitting composition which are defined over fields of characteristic $0$
or $>d$ were classified by Schafer (see [S1] for an overview, [S2] for cubic forms and [S3] for forms of higher degree).
Both the restriction on the dimension of the underlying vector space to be finite, and on
the characteristic of the base field which were needed in Schafer's proofs, were later omitted, respectively relaxed, by McCrimmon [M1].
Forms permitting composition over arbitrary base rings instead of fields were first studied by
Baumgartner and Bergmann [B-B] in 1974, McCrimmon [M2] in 1985 and Petersson [P] in 1993.
 McCrimmon studied quadratic forms permitting composition in the context of nonassociative algebras with a scalar
 involution.
 Petersson generalized the classical Cayley-Dickson doubling process, and completely classified
 composition algebras over arbitrary rings which, as modules, are finitely generated projective with full support.
 Baumgartner and Bergmann  investigated nondegenerate cubic forms over arbitrary commutative unital rings satisfying
a certain multiplicativity condition which is a canonical generalization of what
we call ``permitting composition''. In particular, any unital nonassociative algebra over $R$
carrying a nondegenerate multiplicative cubic form mapping the unit element of the algebra onto the unit of the base ring,
was shown to be alternative and algebraic [B-B, Satz 3, Satz 4].
This result was  extended to arbitrary multiplicative forms of higher degree over rings $R$ with $[d/2]!\in R^\times$
by Legrand-Legrand [L-L].
Since there are no restrictions on the characteristic of the base
ring, the more general concept of a form of degree $d$ as it was developed by Bergmann [B]  was used (see also Roby [R]; the
concept of forms developed by him comprises the one of Bergmann, both coincide for finitely generated
projective modules).
In the setting we consider, their forms are identical with ours and some of Schafer's results carry over verbatim to the setting
of forms permitting composition over rings.

 Petersson was the first to investigate composition algebras over locally ringed spaces, thus initiating the study of nonassociative algebras
 over algebraic varieties. In particular, he generalized the classical Cayley-Dickson doubling
  process due to Albert [A] to this more general setting [P, 2.5]. Composition algebras  (defined on locally free modules of constant finite rank) were
 classified over curves of genus zero in [P, 4.4].

We study  forms of degree $d\geq 3$ permitting composition which are defined over rings with $d!\in R^\times$ instead of fields.
We also investigate unital nonassociative $ \mathcal{O}_X$-algebras $\mathcal{A}$ over locally ringed spaces $X$,
 $d !\in H^0(X,\mathcal{O}_X^\times)$, which carry a nondegenerate form $N:\mathcal{A}\to \mathcal{O}_X$ of degree
 $d\geq 3$ permitting composition,
i.e.,  $N(xy)=N(x)N(y)$ for all sections $x,y$ of $\mathcal{A}$  over the same open subset of $X$
 and $N(1_\mathcal{A})=1$.

The contents of this paper are as follows.
Let $R$ be a commutative ring such that $d!\in R^\times$.
Let $A$ be a unital nonassociative $R$-algebra which is finitely generated projective and faithful as an $R$-module.
Suppose that $A$ carries
 a nondegenerate  form $N$ of degree $d$ permitting composition. Notation and basic facts are given in Section 1.
 After some straightforward generalizations of results in [S3] to forms over $R$ in Section 2,
the cases $d=3$ and $d=4$ are considered separately in Section 4. In particular, we investigate the properties of  $A$
over a domain.
Examples of cubic and quartic forms  over the rings $k[t]$ and $k[t,\frac{1}{t}]$ which permit composition
 are given in Section 5.
In Section 6 our previous results are generalized to the setting of locally ringed spaces.
In Section 7, algebras  over a curve of genus zero which
admit a  nondegenerate cubic form $N$ permitting composition are listed.
As an application, which shows how useful the theory of alternative algebras over varieties can be, we construct classes of unital alternative algebras of degree
$3$ over a field of characteristic not 2 or 3 with a highly degenerate form $N$ which permits
composition in Section 8.

In the following, the standard terminology from algebraic geometry, see Hartshorne's book [H], is used.

\section{Preliminaries}

\subsection{} Let $R$ be a unital commutative associative ring. The {\it rank} of a finitely
 generated projective $R$-module $M$ is defined as ${\rm sup}\{{\rm rank}_{R_P}M_P\,|\,P\in{\rm Spec}\,R\}$.
 Let $A$ be a nonassociative unital $R$-algebra. The {\it nucleus} $\text{ Nuc}(A) = \{ x \in A \, \vert \, [x, A, A] = [A, x, A] = [A, A, x] = 0 \}$
 of  $A$ is the set of those elements in $A$, which associate with all elements in $A$.
The {\it center} of  $A$ is the set of all elements which commute and associate
with all elements of $A$; that is, the set $Z(A)=\{x\in {\rm Nuc}(A)\,|\, xy=yx \text{ for all }y\in A\}$
[S4, p.~14].

 \begin{remark}
 (i) An $R$-module $M$ has {\it full support} if
${\rm Supp}\,(M)=\{P\in {\rm Spec}\,R| M_P\not=0\}={\rm Spec}\,R$.
Every nonassociative unital algebra over $R$ which is finitely generated projective with full support is
{\it faithful}; i.e., $rA=0$ implies $r=0$.
\\(ii) If $A$ is  finitely generated projective and  faithful
as $R$-module, $R1_A$ is a direct summand of $A$:
The map $\epsilon: R\to A$, $\epsilon(r)=r1_A$ has a section if and only if ${\rm Hom}(\epsilon,1_A)$
is surjective.  Since $A_P$ has positive dimension, the map is surjective modulo any maximal ideal, which implies the
assertion (see the proof of [Kn1, I.(1.3.5)]).
\end{remark}

 Unless stated otherwise, the term ``$R$-algebra"  refers to unital nonassociative algebras which are finitely generated
 projective and faithful as $R$-modules.

\subsection{}

 An $R$-algebra $C$ is called a {\it composition algebra}
if it admits a quadratic form $N \colon C \to R$  such that $N(1_C)=1$ which permits composition and
whose induced symmetric bilinear form is {\it  nondegenerate};
i.e., it determines an $R$-module isomorphism $C \overset{\sim}{\longrightarrow}
\check{C} = {\rm Hom}_R (C, R)$.
$C$ is alternative and quadratic.
The quadratic form $N$ is uniquely determined by $C$ and called the {\it norm} of $C$.
 Composition algebras over $R$ only exist in ranks 1, 2, 4 or 8. Those of constant rank~2 are called
{\it quadratic \'etale} $R$-algebras. Those of constant rank 4 are called {\it quaternion algebras} and the ones
of constant rank 8 are called {\it octonion algebras}.
 If there is a  nondegenerate quadratic form $N:A\to R$ on a unital nonassociative $R$-algebra $A$
permitting composition, then $A$ is a composition algebra over $R$ and of rank $1$, $2$, $4$ or $8$ [P].

\ignore{\begin{proof}
see referee's comment, too!!!

 For all $P\in {\rm Spec}\,R$ the residue class algebra $A(P)=A \otimes k(P)$
over the residue class field $k(P)=R_P/m_P$ carries a nondegenerate quadratic form permitting composition and thus is a
composition algebra of dimension $1$, $2$, $4$ or $8$. By [P, 1.7], so is $A$.
\end{proof}}

\subsection{} Let $d$ be a positive integer and assume $d!\in R^\times$. Let $M$ be a finitely generated
projective $R$-module.
A {\it form of degree $d$} (on $M$) over $R$ is a map $N:M\to R$
 such that $N(a x)=a^d N(x)$ for all $a\in R$, $x\in M$, and where the map $\theta : M \times
\dots \times M \to R$ defined by
 $$\theta(x_1,\dots,x_d)=\frac{1}{d!}\sum_{1\leq i_1<\dots<i_l\leq d}(-1)^{d-l}N(x_{i_1}+\dots+x_{i_l})$$
 is a $d$-linear form over $R$ (the range of summation of $l$ being $1\leq l\leq d$). $\theta$ is called the {\it  symmetric $d$-linear form}
associated  with $N$ and $(M,\theta)$ a $d$-linear space.   Obviously, $N(x)=\theta(x,\dots,x)$.
 Note that a symmetric $R$-module homomorphism $M\otimes\dots\otimes M\to R$ or an $R$-module homomorphism
 $S_d(M)\to R$, where $S_d(M)$ is the symmetric algebra of $M$, also define a symmetric $d$-linear form.

 $N$ (or, respectively, the associated $d$-linear form $\theta$) is called {\it  nondegenerate}, if
 the map $ M\to {\rm Hom}_{R}(M\otimes\dots\otimes M,R)$ (($d-1$)-copies of $M$)
$$\theta_{x_1}(x_2\otimes\dots\otimes x_d)=\theta(x_1,x_2,\dots,x_d)$$
is injective (cf. Keet [K, p.~4946]). This concept of nondegeneracy is not invariant under base change.
 A stronger notion, invariant under base change,
is to require that $\theta_{x_1}\otimes k(P)$ is injective for all $P\in {\rm Spec}\,R$. This is equivalent to
saying that $\theta_{x_1}$ is an isomorphism of $M$ onto a direct sumand of ${\rm Hom}_{R}(M^{\otimes (d-1)},R)$.
Both notions, however, are equivalent for  forms permitting composition, see Lemma 1 (i).

Two $d$-linear spaces $(M_i,\theta_i)$, $i=1,2$ are called {\it isomorphic}
(written
$(M_1,\theta_1)\cong (M_2,\theta_2)$ or just $\theta_1\cong\theta_2$)
if there exists an $R$-module isomorphism $f:M_1\to M_2$ such that
$\theta_2(f(x_1),\dots,f(x_d))=\theta_1(x_1,\dots,x_d)$ for all $x_1,\dots,x_d\in
M_1.$

The {\it orthogonal sum} $(M_1,\theta_1)\perp (M_2,\theta_2)$ of $(M_i,\theta_i)$, $i=1,2$, is defined to be the $R$-module
 $M_1\oplus M_2$ together with the $d$-linear form
$(\theta_1 \perp \theta_2)(u_1+x_1,\dots,u_d+x_d)=\theta_1(u_1,\dots,u_d)+\theta_2(x_1,\dots,x_d)$.
A $d$-linear space $(M,\theta)$ is
called {\it decomposable}, if $(M,\theta)\cong (M_1,\theta_1)\perp (M_2,\theta_2)$
for two non-zero $d$-linear spaces $(M_i,\theta_i)$, $i=1,2$.
A non-zero $d$-linear space $(M,\theta)$ is
called {\it indecomposable} if it is not decomposable. We  distinguish between indecomposable ones and
{\it absolutely} indecomposable ones;
i.e., $d$-linear spaces which stay indecomposable under base change.

\subsection{} Let $X$ be a locally ringed space with structure sheaf
$\mathcal{O}_X$. For $P \in X$ let $\mathcal{O}_{P,X}$ be the local ring of
$\mathcal{O}_X$ at $P$ and $m_P$ the maximal ideal of $\mathcal{O}_{P,X}$. The corresponding residue class
field is denoted by $k(P)=\mathcal{O}_{P,X}/m_P$. For an $\mathcal{O}_X$-module $\mathcal{F}$ the stalk of
$\mathcal{F}$ at $P$ is denoted by $\mathcal{F}_P$. $\mathcal{F}$ is said to have {\it full support} if
${\rm Supp}\,\mathcal{F}=X$; i.e., if $\mathcal{F}_P\not=0$ for all $P\in X$. We call
$\mathcal{F}$ {\it locally free of finite rank} if for each $P\in X$ there is an open neighborhood $U\subset X$
of $P$ such that $\mathcal{F}|_U=\mathcal{O}_U^r$ for some integer $r\geq 0$. The {\it rank} of
$\mathcal{F}$ is defined to be ${\rm sup}\{{\rm rank}_{\mathcal{O}_{P,X}}\mathcal{F}_P\,|\, P\in X\}$.
 The term ``$\mathcal{O}_X$-algebra" (or ``algebra over $X$'') always refers to nonassociative
$\mathcal{O}_X$-algebras which are unital and locally free of finite rank as $\mathcal{O}_X$-modules.

\subsection{}
Let $\mathcal{C}$ be an $\mathcal{O}_X$-algebra. $\mathcal{C}$ is called a {\it composition algebra} over $X$
if it has full support and if there exists a {\it nondegenerate} quadratic form $N \colon \mathcal{C} \to  \mathcal{O}_X$
(i.e., the induced symmetric bilinear form $N(u,v)  =
N(u+v)-N(u)-N(v)$ determines a module isomorphism
$\mathcal{C} \overset{\sim}{\longrightarrow} \check{\mathcal{C}}
=\mathcal{H}om (\mathcal{C},\mathcal{O}_X)$), such that $N(uv)=N(u)N(v)$ for all sections $u,v$ of $\mathcal{C}$ over the
same open subset of $X$ [P, 1.6].

The form $N$ is uniquely determined by these conditions and called the {\it norm} of $\mathcal{C}$. It is denoted by $N_\mathcal{C}$.
Given an algebra $\mathcal{C}$ over $X$ and a quadratic form $N \colon \mathcal{C}
\to \mathcal{O}_X$, the algebra $\mathcal{C}$ is a composition algebra over $X$ with
norm $N$ if and only if $\mathcal{C}_P$ is a composition
algebra over $\mathcal{O}_{P,X}$ with norm $N_P$ for all $P \in X$.
 Composition algebras over $X$ are invariant under base change, and exist only in ranks 1, 2, 4 or 8.
A composition algebra of constant rank 2 (resp. 4 or 8) is called a {\it quadratic \'etale algebra} (resp. {\it quaternion algebra} or
an {\it octonion algebra}). A composition algebra over $X$ of constant rank is called {\it split}, if it contains a composition subalgebra
isomorphic to $\mathcal{O}_X \oplus \mathcal{O}_X$ [P, 1.7, 1.8].

 If $X$ is an $R$-scheme with structure morphism
$\tau \colon X \to {\rm Spec}\,R$, then a composition algebra $\mathcal{C}$ over $X$
is  {\it defined over} $R$ if there exists a composition algebra $C$ over $R$ such that
$\mathcal{C}=\tau^{\ast} C \cong C \otimes \mathcal{O}_X$. There exists a generalized
Cayley-Dickson doubling ${\rm Cay} (\mathcal{D},\mathcal{P},N)$
 for a composition algebra $\mathcal{D}$ of constant rank $\leq 4$ over a locally ringed space [P, 2.3, 2.4, 2.5].

\subsection{} Let $d !\in H^0(X,\mathcal{O}_X^\times)$. Let $\mathcal{M}$ be an $\mathcal{O}_X$-module which is locally free of finite rank.
A {\it form of degree $d$} (on $\mathcal{M}$) over $\mathcal{O}_X$ is a map $N:\mathcal{M}\to \mathcal{O}_X$
 such that $N(a x)=a^d N(x)$ for all sections $a$ of $\mathcal{O}_X$, $x$ of $\mathcal{M}$ over the same open subset
 of $X$, and where the map $\theta : \mathcal{M} \times \dots \times \mathcal{M}\to \mathcal{O}_X$ defined by
 $$\theta(x_1,\dots,x_d)=\frac{1}{d!}\sum_{1\leq i_1<\dots<i_l\leq d}(-1)^{d-l}N(x_{i_1}+\dots+x_{i_l})$$
($1\leq l\leq d$) is a $d$-linear form over $\mathcal{O}_X$. $\theta$ is called the {\it  symmetric $d$-linear form}
associated  with $N$. (This definition of associated $d$-linear form deviates from the one used for
the associated symmetric bilinear form in 2.5 or [P] by the factor $1/2$ which was omitted in order to be able to include the
 case that $2\not\in H^0(X,\mathcal{O}_X^\times)$ in the classification result for composition algebras).
A form $N$ of degree $d\geq 3$ (or, respectively, its associated $d$-linear form $\theta$) is called
{\it  nondegenerate}, if
 the map $ \mathcal{M}\to {\mathcal H}om_X(\mathcal{M}\otimes\dots\otimes \mathcal{M},\mathcal{O}_X)$ ($(d-1)$-copies of
 $\mathcal{M}$),
$$\theta_{x_1}(x_2\otimes\dots\otimes x_d)=\theta(x_1,x_2,\dots,x_d)$$
is injective.

Two $d$-linear spaces $(\mathcal{M}_i,\theta_i)$ ($i=1,2$) are called {\it isomorphic}
if there exists an $\mathcal{O}_X $-module isomorphism $f:\mathcal{M}_1\to \mathcal{M}_2$ such that
$\theta_2(f(v_1),\dots,f(v_d))=\theta_1(v_1,\dots,v_d)$ for all sections $v_1,\dots,v_d$ of
$\mathcal{M}_1$ over the same open subset of $X$.

The {\it orthogonal sum} $(\mathcal{M}_1,\theta_1)\perp (\mathcal{M}_1,\theta_2)$ of $(\mathcal{M}_i,\theta_i)$, $i=1,2$, is defined to be the $\mathcal{O}_X$-module
 $\mathcal{M}_1\oplus \mathcal{M}_2$ together with the $d$-linear form
$(\theta_1 \perp \theta_2)(u_1+v_1,\dots,u_d+v_d)=\theta_1(u_1,\dots,u_d)+\theta_2(v_1,\dots,v_d)$.
A $d$-linear space $(\mathcal{M},\theta)$ is
called {\it decomposable}, if $(\mathcal{M},\theta)\cong (\mathcal{M}_1,\theta_1)\perp (\mathcal{M}_2,\theta_2)$
for two non-zero $d$-linear spaces $(\mathcal{M}_i,\theta_i)$, $i=1,2$.
A non-zero $d$-linear space $(\mathcal{M},\theta)$ is
called {\it indecomposable} if it is not decomposable. We will distinguish between indecomposable ones and {\it absolutely}
 indecomposable ones;
i.e., $d$-linear spaces which stay indecomposable under base change.

\section{Forms permitting composition}

Large parts of the results and proofs of [S3] can be generalized verbatim to the case where the base field
is replaced by an arbitrary
commutative ring $R$ with $d!\in R^\times$. We briefly summarize them for the sake of the reader.

Let $R$ be a ring with $d!\in R^\times$.
Unless explicitly stated otherwise, the term ``$R$-algebra" refers to unital nonassociative
algebras which are finitely generated projective and faithful as $R$-modules.
Let $A$ be an  algebra over $R$,  $1=1_A$ the unit element of $A$ and $N:A\to R$ a form of degree $d$ permitting composition; i.e.,
$$\theta(xy,\dots,xy)=\theta(x,\dots,x)\theta(y,\dots,y)$$
for all $x,y\in A$ and $N(1_A)=1$.
 Linearizing this equation in $x$ and then in $y$ we obtain
$$\theta(x_1y,\dots,x_dy)=\theta(x_1,\dots,x_d)N(y),$$
and
$$\sum_\sigma \theta(x_1y_{\sigma(1)},\dots,x_d y_{\sigma(d)})=d! \theta (x_1,\dots,x_d)\theta(y_1,\dots,y_d),$$
where $\sigma$ ranges over all the permutations in $S_d$. This implies
$$\theta(xy_1,\dots,xy_d)=N(x)\theta(y_1,\dots,y_d)$$
by symmetry [S3, (5), (6) and (7)]. (Indeed, this observation does not require $A$ to have a unit element.)
 For $i=1,\dots,d$ we define a form $T_i:A\to R$ of degree $i$ via
$$T_i(x)=\binom{d}{i}\theta(x,\dots, x,1,\dots,1)\,\,\text{  ($i$-times }x).$$ Then
$$N(x)=T_d(x) \text{ and } T_1(x)=d\theta(x,1,\dots,1).$$
Put $T_0(x)=1$ and $T_{d+q}(x)=0$ for $q>0$.
The form $T:A\to R$, $T(x)= T_1(x)$ is called the {\it trace}. Define $A_0={\rm ker}\, T$. $A$ can be written as the direct sum of $R$-modules
$$A=R1_A\oplus A_0.$$

From now let $d\geq 3$, unless explicitly stated otherwise.

A $d$-linear form $\theta$ is {\it invariant under
all left and right multiplications with elements of trace zero},
if
$$\theta(x_1a,\dots,x_d)+\theta(x_1,x_2a,\dots,x_d)+\cdots +\theta(x_1,\dots,x_da)=0$$
and
$$\theta(ax_1,\dots,x_d)+\theta(x_1,ax_2,\dots,x_d)+\cdots +\theta(x_1,\dots,ax_d)=0$$
for all $x_i\in A$ and for all elements $a\in A$ of trace zero.

 $B(x,y)=T(xy)$ is a symmetric  bilinear form on $A$ which is associative, and
 if $N$ is nondegenerate, then so is $B$.
 $A$ is alternative; i.e., $x^2y=x(xy)$ and $yx^2=(yx)x$ for all elements $x,y\in A$. Every element $x\in A$ satisfies
$$x^d-T_1(x)x^{d-1}+T_2(x)x^{d-2}-\cdots+(-1)^d T_d(x)1=0$$
([S3, Theorem 2, Theorem 3], for $d=2$ this was shown in [M2, 4.6]).

\begin{remark} (i) If  $N$ is degenerate, the algebra $A$ need not be alternative; for a counterexample for $d=2$ see [M2, 4.13].
\\(ii) Baumgartner and Bergmann studied cubic forms $N$ on unital nonassociative
algebras $A$ over arbitrary
rings which they called multiplicative, and called such a pair $(A,N)$ a {\it composition algebra of third degree}.
Each nondegenerate composition algebra (in the sense of [B-B]) of degree 3 over an arbitrary ring is alternative
 [B-B, Satz 3].  This result for composition algebras of third degree
 in the sense of [B-B] was later generalized by Legrand-Legrand [L-L]
 to multiplicative forms of higher degree over rings with $[d/2]!$ invertible.
 Their definition of forms is different from ours in order to accommodate the case that $d!$ is not
 invertible in $R$.
 Over fields of characteristic not 2 or 3, or over the rings  considered here,
  the different concepts of forms coincide and a form $N$ permitting composition which satisfies $N(1_A)=1$
 corresponds to a multiplicative form.  For the above equation in a more general context, see [L-L, Corollaire 2.2].
\end{remark}

\begin{theorem}  Suppose that $R$ is a domain. Let $A$ be an $R$-algebra,
 $\theta:A\times\dots\times A\to R$ a nondegenerate symmetric $d$-linear form on $A$
and let the trace $T:A\to R$ be defined as above. If $\theta$ is invariant under all left and right multiplications
corresponding to elements of trace zero, then $A$ is semiprime; i.e.,  $I^2\not=0$
for each non-zero ideal $I$ in $A$.
\end{theorem}

For $d>2$ the proof can be found in the proof of  [S3, Theorem 1 (c)]. A similar argument  holds for $d=2$.

\begin{corollary} Suppose that $R$ is a domain and $A$ an  $R$-algebra with  be a nondegenerate
form $N:A\to R$ of degree $d\geq 2$ permitting composition. Then $A$ is a semiprime alternative algebra over $R$.
If $I$ is a minimal ideal of $A$ and $I^2\not=(0)$, then $I$ is simple and $I^2=I$. Moreover, $I$ is either a
simple associative ring or an octonion algebra over its center. In the latter case, $A=I\oplus I'$ holds as an
ideal direct sum.
\end{corollary}

\begin{proof}
 The first assertion follows immediately from Theorem 1 using that the $d$-linear form associated with $N$ is
 invariant under all left and right multiplications
corresponding to elements of trace zero, the second from [Sl2, Theorem A, Lemma 2.1, Theorem B].
\end{proof}

\begin{remark} Suppose $d\geq 2$ and that the ring  $A$ satisfies the descending chain condition on right ideals.
 By Corollary 1, $A$ is a semiprime alternative algebra over $R$. By [Sl1, Theorem B], $A$ is expressible,
 unique up to order, as an ideal direct sum $A=A_1\oplus \cdots \oplus A_r$
of minimal ideals $A_i$ where each $A_i$ is either an octonion algebra over its center
or a simple artinian associative ring. If every non-zero ideal of $A$ contains a minimal ideal of $A$ and $A$ is purely alternative (i.e.,
has no non-zero nuclear ideals), then $A$ is expressible, unique up to order, as an ideal direct sum
$A=A_1\oplus \cdots \oplus A_r$ where each $A_i$ is an octonion algebra over its center [Sl2, Theorem A].
\end{remark}

\begin{lemma} (i) $N(P)$ is nondegenerate for all $P\in {\rm Spec}\, R$.
\\ (ii) If $A$ is an associative algebra
then $A$ is separable. If $A$ is also central then it is an Azumaya algebra over $R$.
\end{lemma}

\begin{proof} (i) Since $N$ is nondegenerate, so is $T$ and thus also $T(P)$ for all $P\in X$.
 The fact that $T(P)$ is nondegenerate is equivalent to $N(P)$ being nondegenerate ([B, Satz 3] or [B-B, p.~327] for cubic
 forms) implying
 that $N(P)$ is a nondegenerate  form for all $P\in X$.
 \\(ii) $A$ is finitely generated as $R$-module by assumption, thus $A$ is separable over $R$ if and only if $A/mA$
is a separable $R/m$-algebra for all $m\in {\rm Max}\,R$ [Kn1, III. (5.1.10)].
This holds, since  $A(P)$ admits a nondegenerate form of degree $d$ permitting composition for all
$P\in {\rm Spec}\, R$. Thus $A(P)$ is a finite dimensional separable $k(P)$-algebra
for all $P\in {\rm Spec}\, R$ by [S1, Theorem 3]. Hence $A$ is a separable associative $R$-algebra,
finitely generated as an $R$-module. If $A$ is central, then $A$ is central, separable and finitely generated
 faithfully projective as an $R$-module, therefore an Azumaya algebra by [Kn1, III. (5.1.1)].
\end{proof}

Thus, if $N$ permits composition, by (i) the two different notions of nondegeneracy introduced in 2.3 are equivalent.

\begin{remark} Assume $A$ has  center $R'$ larger than $R$ in the situation of Lemma 1 (ii).
Then $R'$ is a separable ring extension of $R$.
 View $A$ as an $R'$-algebra. $A$ is finitely generated and faithful as an $R'$-module.
 In case $A$ is projective  and separable also as an $R'$-module, $A$ is an Azumaya algebra over $R'$
 [Kn1, III, (5.1.1)].
\end{remark}

\begin{proposition} $A(P)$ is a finite-dimensional separable algebra over $k(P)$ and\\
(i) if $N$ is a cubic form, then $A$ has rank $1$, $2$, $3$, $5$ or $9$ over $R$;\\
(ii) if $N$ is a quartic form, then $A$ has rank $1$, $2$, $3$, $4$, $5$, $6$, $8$, $9$, $10$, $12$ or $16$ over $R$.
\end{proposition}
\begin{proof} For all $P\in {\rm Spec}\,R$ the residue class algebra $A(P)$ is  a non-zero alternative algebra
 over the residue class field $k(P)$ together with a nondegenerate form $N(P):A(P)\to k(P)$ which permits composition.
Therefore $A(P)$ is a finite-dimensional separable algebra over $k(P)$ [S1, Theorem 3].
 Hence, if $N$ is a cubic form, then $A$ must have rank
$1$, $2$, $3$, $5$ or $9$, and if it is a quartic form then $A$ must have rank
$1$, $2$, $3$, $4$, $5$, $6$, $8$, $9$, $10$, $12$ or $16$ [S2, S3].
\end{proof}

Corresponding statements  can be derived for forms of higher degree than 4.

\subsection{} For this section we assume that $R$ is a domain.

\begin{lemma} Suppose that $A$ contains an idempotent $e\not=0,1$. Then $1$ and $e$ are linearly
independent over $R$.
\end{lemma}

\begin{proof} Let $1a+eb=0$ for $a,b\in R$. Multiplication with $e$ implies  $(a+b)e=0$ and hence $a+b=0$, since
$A$ is a projective $R$-module, thus torsion free as an $R$-module. Now $b=-a$, so $0=1a+eb=(1-e)a$ yields
$a=0$ since $1-e\not=0$, hence also $b=0$.
\end{proof}

 Suppose that $A$ contains an idempotent $e\not=0,1$. Then
 $$x^d-T_1(x)x^{d-1}+T_2(x)x^{d-2}-\cdots+(-1)^d T_d(x)1=0$$ for all $x\in A$  implies that
 $$1-T(e)+T_2(e)-\cdots+(-1)^{d-1} T_{d-1}(e)1=0,$$
and $N(e)=0$, since $1$ and $e$ are linearly independent. As in [S3, p.~785] we obtain the equations
 $$(j+1)T_{j+1}(e)=(T(e)-j)T_j(e)$$
 for $j=1,\dots,d-1$ [S3, (38)], and
 $$T_{i}(e)=\binom{m}{i}$$
 for $i=1,\dots,m$ [S3, (40)],
 where $m$ is the least integer such that $T_{m+1}(e)=\cdots=T_d(e)=0.$

 These are needed for the proof of the next result, which closely follows the one given in [S3] for Theorem 4.

 \begin{theorem}  If
$$A=A_1\oplus \cdots \oplus A_r$$
is the direct sum of ideals $A_i\not=0$ in $A$, then
 $$N(x)=N_1(x_1)\cdots N_2(x_r)$$
 where $x=x_1+\dots +x_r$, $x_i\in A_i$. Each $N_i$ is a nondegenerate form of degree $d_i$ on $A_i$
 which permits composition and $d=d_1+\dots+d_r$. If $r\geq 2$ then $N$ is absolutely indecomposable.
 \end{theorem}

 \begin{proof} Assume that $A=G\oplus G'$ with $G\not=0$, $G'\not=0$ ideals.
 Write $1=e+e'$ with $e\in G$ and $e'\in G'$. Then $e\not=0,1$ (resp. $e'\not=0,1$) is the unit element
 of $G$ (resp. of $G'$) and as such it is an idempotent in $A$. For $g\in G$, $g'\in G'$ define
 $$N_G(g)=N(g+e),\,\,N_{G'}(g')=N(e+g').$$
 For any $x=g+g'\in A$ we get
 $$N(x)=N_G(g)N_{G'}(g')$$
 and $$N_G(g_1g_2)=N_G(g_1)N_G(g_2)$$
 for all $g_1,g_2\in G$.
 By showing that
 $$N_G(g)=T_m(g) \text{ for all }g\in G,\text{ where }m=T(e),$$
 $N_G$ is proved to be a form of degree $m$ over $G$.
 (The fact that $1$ and $e$ are linearly independent is needed to obtain $T(e)=m$.)
 Symmetrically, the same formulas hold for $N_{G'}$, so that $N(x)$ is a product of forms of degree $m$ and $m'$
  permitting composition and $d=m+m'$.
  The proof that $N_G(g)=T_m(g) \text{ for all }g\in G,\text{ with }m=T(e)$ is the same as given in [S3].
Both $N_G$ and $N_{G'}$ are nondegenerate, see the proof in [S3, p.~787, 788].

This argument can be repeated finitely often and since $A=A_1\oplus \cdots \oplus A_r$ with $A_i\not=0$
 ideals, we obtain the assertion by induction.\\
 That $N$ is absolutely indecomposable follows from [Pr, 5.1].
 \end{proof}

If $R$ is a field, it is well-known that
 any form  on a simple alternative algebra which permits composition is a power
 of the generic norm of the algebra. If $R$ is a ring, it can happen that an algebra admits
 more than one nondegenerate form permitting composition, see [M2].

 \begin{remark} (i) Let $A=A_1\oplus A_2$ be a direct sum of two ideals.
 If $A$ is finitely generated projective as an $R$-module then so are $A_1$ and $A_2$.
 \\ (ii) If $d=3$ (resp. $4$) in Theorem $2$, then $A$
 is the direct sum of at most three (resp. four) non-zero ideals which all are unital nonassociative algebras
  admitting a form  of degree
 $1$ or $2$ (resp. $1$, $2$ or $3$) permitting composition.
\end{remark}

\section{Cubic and quartic forms}

 A nondegenerate cubic (or quartic) form $N:A\to R$ permitting composition
can only exist in rank $1$, $2$, $3$, $5$ or $9$ (or in rank
$1$, $2$, $3$, $4$, $5$, $6$, $8$, $9$, $10$, $12$ or $16$) by Proposition 1.
 [S3, Lemma 1, 2, 3] applied to the residue class algebras imply the next two lemmas:

\begin{lemma}  Suppose that there exists a cubic form $N:A\to R$
 permitting composition.\\
(i) If $A$ is a non-split quadratic \'{e}tale algebra over $R$, then  the residue class algebras $A(P)$
must be split for all $P\in{\rm Spec}\,R$.\\
(ii)  $A$ cannot be a quaternion or octonion algebra over $R$.
\end{lemma}

\begin{lemma}  Let $A$ be a composition algebra over $R$ of constant rank  with norm $n_A$.
Suppose that there exists a quartic form $N:A\to R$
 permitting composition.\\ 
 (i) If $A$ is a  quadratic \'{e}tale algebra and $A(P)$ is a quadratic field extension for all $P\in{\rm Spec}\,R$, or
  if $A$ has rank $\geq 4$  then $N(P)(x)=n_{A(P)}(x)^2$  for all $P\in{\rm Spec}\,R$ and $x\in A(P)$.\\
  (ii) $A$ cannot be an Azumaya algebra over $R$ of rank 9.\\
  (iii) If $A$ is a cubic ring extension  of $R$, then $A(P)$ is not a cubic field extension of $k(P)$.
  for all $P\in{\rm Spec}\,R$.
\end{lemma}

\begin{theorem} (i) Let $A$ be an $R$-algebra of constant rank
  such that there exists a  cubic form $N$ on $A$ permitting composition. Suppose that for each
$P\in {\rm Spec}\,R$ there exists an element $u\in A\otimes_R k(P)$ such that $1,u,u^2$ are linearly
independent over $k(P)$ (by Theorem 2, $A$ is alternative, so the powers of $A$ are unambiguous).  Let $M$, $Q$ and $L$ be a cubic, a quadratic and a linear form from $A$ to $R$
satisfying $$x^3-L(x)x^{2}+Q(x)x^{}-M(x)1=0$$ for all $x\in A$.
Then $M=N$, $S=T_2$ and $L=T$.\\
(ii) Let $A$ be an $R$-algebra of constant rank such
that there exists a quartic form $N$ on $A$ permitting composition. Suppose that for each
$P\in {\rm Spec}\,R$ there exists an element $u\in A\otimes_R k(P)$ such that $1,u,u^2,u^3$ are linearly
independent over $k(P)$.  Let $S$, $M$, $Q$ and $L$ be a quartic, cubic, quadratic and a linear form from $A$ to $R$
satisfying $$x^4-L(x)x^{3}+Q(x)x^{2}-M(x)x+S(x)1=0$$ for all $x\in A$.
Then $S=N$, $M=T_3$, $Q=T_2$ and $L=T$.
\end{theorem}

Part (i) is a Corollary of [Ach, 1.12] applied to the Jordan algebra $A^+$ determined by $A$, (ii) can be proved analogously (if more tediously). Indeed, both results remain true
even if we remove the restriction on $R$ to satisfy $d!\in R^\times$ and work with the more general notion of a form of
higher degree as given in [R].

\begin{remark} Suppose that $A$ is an Azumaya algebra of constant rank $9$ over $R$. For each
$P\in {\rm Spec}\,R$ there exists an element $u\in A(P)$ such that $1,u,u^2$ are linearly
independent over $k(P)$, since  $A(P)$ is a central simple algebra of degree 3.
 Since the reduced norm $n$ and trace $t$ of $A$ satisfy
$x^3-t(x)x^{2}+q(x)x^{}-n(x)1=0$ with $q$ a quadratic form, it follows that $M=n$ (and $L=t$).

Suppose that $A$ is a cubic  ring extension of $R$.
 For each $P\in {\rm Spec}\,R$ there exists an element $u\in A(P)$ such that $1,u,u^2$ are linearly
independent over $k(P)$, since $A(P)$ is a cubic  \'etale algebra over $k(P)$.
 Since the reduced norm $n$ and trace $t$ of $A$ satisfy
$x^3-t(x)x^{2}+q(x)x^{}-n(x)1=0$  with $q$ a quadratic form, it follows again that $M=n$ (and $L=t$).

Analogous arguments show that also for a quartic separable ring extension of $R$ and for an
 Azumaya algebra of constant rank $16$ over $R$, any quartic form on $A$ which permits composition must be uniquely determined
 and be the norm of the algebra, if $k(P)$ is an infinite field for all $P\in {\rm Spec}\,R$.
\end{remark}

For the rest of this section, let $R$ be a domain.

\begin{proposition} Let $N:A\to R$ be a cubic form on $A$ permitting composition.
Then $(A,N)$ is one of the following:\\
 (i) $A=R$ and $N(x)=x^3$;\\
 (ii) $A$ is a separable commutative associative $R$-algebra of rank 2 or 3, and $N$ is absolutely
 indecomposable; if $A$ is a non-split quadratic \'{e}tale algebra over $R$, then the residue class algebras $A(P)$
must be split for all $P\in{\rm Spec}\,R$, if $A$ has rank 3, then $A(P)$ is a  cubic  \'{e}tale algebra over $k(P)$
 for all $P\in {\rm Spec}\,R$.\\
(iii) $A$ has rank $5$ and is a separable associative, but not commutative, $R$-algebra and $N$ is absolutely indecomposable.\\
(iv) If  $A$ has rank $9$ and is associative, then it is an Azumaya algebra over $R$ and $N$ is its
- uniquely determined - norm.
If $A$ has rank $9$ and is not associative, then $A$ is not commutative and $N$ is absolutely indecomposable.
\end{proposition}

\begin{proof}
 Let $m$ denote the rank of $A$, then $m\in\{1,2,3,5,9\}$.
\begin{enumerate}
\item{} If $m=1$, then $A=R$ and $N(x)=x^3$.

\item{} If $m=2$, then $A(P)\cong k(P)\oplus k(P)$ is commutative, associative  and $N(P)(x_1+x_2)=x_1x_2^2$ is absolutely
indecomposable for all $P\in {\rm Spec}\,R$. Thus all commutators and associators lie in $IA$ where $I$ is the nil radical
of $R$. Since $R$ is a domain, $I=0$ and $A$ itself must be a commutative associative $R$-algebra.
 $N(P)(x_1+x_2)=x_1x_2^2$ is absolutely
indecomposable for all $P\in {\rm Spec}\,R$, hence $N$ is an absolutely indecomposable form. The rest follows from  Lemma 3.

\item{} If $m=3$, then $A(P)$ is a  cubic  \'{e}tale algebra over $k(P)$  and $N(P)$ is its
(absolutely indecomposable) norm  for all $P\in {\rm Spec}\,R$.
Therefore $N$ is absolutely indecomposable
 and, by the same argument as above, $A$ must be commutative and associative.

\item{} If $m=5$, then
$A(P)\cong k(P)\oplus \text{ ``some quaternion algebra over } k(P)\text{''}$
for all $P\in {\rm Spec}\,R$. Therefore $A(P)$ is  associative and not commutative and
$N(P)$ is absolutely indecomposable  for all $P\in {\rm Spec}\,R$.
 Thus $A$ is  associative by the same argument as above, and not commutative.  $N$ is absolutely indecomposable.

\item{} If $m=9$ we distinguish two cases:
If $A$ is associative, then $A(P)$ is a central simple algebra of degree $3$ over $k(P)$ for all $P\in {\rm Spec}\,R$
and $N(P)$ is its norm. Therefore $A$ is an Azumaya algebra over $R$ and $N$ is its - uniquely determined - norm.\\
If $A$ is not associative, then
$A(P_0)\cong k(P_0)\oplus \text{ ``some octonion algebra over } k(P_0)\text{''}$
for some $P_0\in {\rm Spec}\,R$ and
$N(P_0)$ is absolutely indecomposable. Thus $A$ is not associative and not commutative
and $N$ is again absolutely indecomposable.
\end{enumerate}
\end{proof}

 \begin{example} Over rings, the first Tits construction
[P-R, Theorem 3.5] starting with $R$ can be generalized as follows [Ach, 2.25]:
Let $L\in \, _3 {\rm Pic} R$ and $N_L:L\to R$ a nondegenerate cubic form, let $L^\vee=
{\rm Hom}_R(L,R)$ be its dual and $\langle w,\check{w}\rangle=\check{w}(w)$ the canonical pairing $L\times\check{L}\to
R.$ There exists a uniquely determined cubic norm $\check{N}_L:L^\vee\to R$ and
uniquely determined adjoints $\sharp:L\to L^\vee$ and $\check{\sharp}:L^\vee\to L$ such that
\begin{enumerate}
\item $\langle w,w^\sharp\rangle=N_L(w)1$;
\item $\langle \check{w}^{\check{\sharp}},\check{w}\rangle=\check{N_L}(\check{w})1$;
\item $w^{\sharp \, \check{\sharp}}=N_L(w)w$;
\end{enumerate}
for $w$ in $L$, $\check{w}$ in $L^\vee$. Moreover,
\begin{enumerate}
\item $\check{w}^{\check{\sharp}\,\sharp }=\check{N_L}(\check{w})\check{w}$;
\item $\langle w, \check{w}\rangle^2= \langle \check{w}^{\check{\sharp}}, w^{\sharp}\rangle$;
\item $\langle w, \check{w}\rangle^3=N_L(w)\check{N}_L(\check{w})$;
\item $\langle w,\check{w}\rangle w= 3\langle w,\check{w}\rangle w-w^\sharp\check{\times}\check{w}$
\end{enumerate}
for $w,w'$ in $L$, $\check{w}$ in $L^\vee$ [Ach, 2.13]. Define
$$\begin{array}{l}
\widetilde{\J}=R\oplus L\oplus L^\vee,\\
\widetilde{1}=(1,0,0),\\
\widetilde{N}(a,w,\check{w})=a^3+N_L(w)+\check{N_L}(\check{w})-3a \langle w,\check{w}\rangle\\
(a,w,\check{w})^{\widetilde{\sharp}}=(a^2-\langle w,\check{w}\rangle, \check{w}^{\check{\sharp}}-
aw, w^\sharp-\check{w}a)
\end{array}$$
for $a\in A$, $w\in L$, $\check{w}\in L^\vee$, then $(\widetilde{N},\widetilde{\sharp},\widetilde{1})$
is a cubic form with adjoint and base point on $\widetilde{\J}$ and has trace form
$$\widetilde{T}((a,w,\check{w}),(c,v,\check{v}))=3ac+3\langle w, \check{v}\rangle+3\langle v,
\check{w}\rangle.$$
The Jordan algebra $A=
\J(\widetilde{N},\widetilde{\sharp},\widetilde{1})$ over $R$
 is commutative and associative
  and the nondegenerate cubic form $\widetilde{N}$ permits composition.
 $A(P)$ is
 a  cubic  \'{e}tale algebra over $k(P)$ for all $P\in {\rm Spec}\,R$.
\end{example}

\begin{lemma} Let $N:A\to R$ be a nondegenerate cubic form on $A$ permitting composition.\\
 (i) Suppose $A$ can be written as the direct sum of two non-zero ideals $A_1$, $A_2$ of $A$.
 Then $A=R\oplus A_2$ where $A_2$ is a composition algebra over $R$ with norm $n$
   and $N(x_1+x_2)=x_1n(x_2)$ is absolutely indecomposable.
\\(ii)  Suppose $A$ can be written as the direct sum of three non-zero ideals $A_1$, $A_2$, $A_3$ of $A$. Then
$A=R\oplus R\oplus R$ and $N(x_1+x_2+x_3)=x_1x_2x_3$.
\end{lemma}

\begin{proof}\begin{enumerate}
\item{} Suppose $A$ can be written as the direct sum of two non-zero ideals.
 Then $A=A_1\oplus A_2$ and $N(x_1+x_2)=n_1(x_1)n(x_2)$ with $n_1(x_1)=x_1$ and $n$
 a nondegenerate quadratic form permitting composition. Therefore
$A=R\oplus A_2$ where $A_2$ is an algebra over $R$ with a nondegenerate quadratic form $n$  permitting composition and
$N(x_1+x_2)=x_1n(x_2)$ is absolutely indecomposable. Since $R$ is a domain,
 $A_2$ has full support, hence it is a composition algebra over $R$.
\item{}  Suppose $A$ can be written as the direct sum of three non-zero ideals $A=A_1\oplus A_2\oplus A_3$. Then
$N(x_1+x_2+x_3)=n_1(x_1)n_2(x_2)n_3(x_3)$ with $n_i$ linear forms permitting composition, $i=1,2,3$. Hence
$A=R\oplus R\oplus R$ and $N(x_1+x_2+x_3)=x_1x_2x_3$.
\end{enumerate}
\end{proof}

From now on suppose that  $k(P)$ is an infinite field for all $P\in {\rm Spec}\,R$.

\begin{proposition} Let $N:A\to R$ be a quartic form on $A$ permitting composition.
Then one of the following holds:\\
 (i) $A=R$ and $N(x)=x^4$.\\
(ii) $A$ is a commutative associative separable $R$-algebra of rank 2 or 3.
If $A$ has rank 3, then $A(P)$ cannnot be a cubic field extension of $k(P)$
  for all $P\in{\rm Spec}\,R$.\\
(iii)  $A$ has rank $4$ and is an associative separable $R$-algebra.
If, in particular, $A$ is not commutative, then  $A(P_0)$ is a quaternion algebra over $k(P_0)$ for at least one $P_0$.\\
(iv)  $A$ has rank $5$ or $6$ and $A$ is an associative separable $R$-algebra, but not commutative.\\
(v) $A$ has rank $8$ and $A$ is not commutative.
If $A$ is not associative, then $A(P_0)$ is an octonion algebra over $k(P_0)$ for some $P_0\in {\rm Spec}\,R$.\\
(vi) $A$ has rank $9$, $10$ or $12$,  is not associative, not commutative.\\
(vii) $A$ has rank $16$ and $A$ is not associative. Then
$A(P_0)$ is an octonion algebra over some quadratic field extension of $k(P_0)$, for some $P_0\in {\rm Spec}\,R$.\\
(viii)  $A$ has rank $16$ and $A$ is an Azumaya algebra over $R$ with norm $N$.\\
In cases (i) to (vii), $N$ is absolutely indecomposable.
\end{proposition}

\begin{proof} Let $m$ be the rank of $A$, then $m\in \{1,2,3,4,5,6,8,9, 10,12,16\}$.

\begin{enumerate}
\item{} If $m=1$, then $A=R$ and $N(x)=x^4$.
\item{} If $m=2$, then $A(P)$ is isomorphic to a quadratic \'{e}tale algebra over $k(P)$
and $N(P)=(x_1+x_2)=x_1^2x_2^2$ is absolutely
indecomposable for all $P\in {\rm Spec}\,R$, implying that $N$ must be absolutely indecomposable.
 By the same argument as used in the proof of Proposition 2 (2), $A$ is a commutative
associative $R$-algebra together.
\item{} If $m=3$, then
$A(P)\cong k(P)\oplus \text{ ``some quadratic \'{e}tale algebra over } k(P)\text{''}$
by [S3, Lemma 2] and $N(P)$ is absolutely indecomposable for all $P\in {\rm Spec}\,R$.
Thus $A$ is commutative,  associative and separable and $N$ absolutely indecomposable.
\item{} If $m=4$, then $A(P)$ is either a quaternion algebra over $k(P)$, a separable quartic field extension
over $k(P)$, a quadratic \'{e}tale algebra over some quadratic field extension of $k(P)$ or
the direct sum of two quadratic \'{e}tale algebras over $k(P)$. Hence $A$ is associative and separable.
 In particular, if $A$ is not commutative, then $A(P_0)$ must be a quaternion algebra over $R$ for a $P_0\in {\rm Spec}\,R$.
\item{} If $m=5$, then $A(P)\cong k(P)\oplus \text{ ``some quaternion algebra over } k(P)\text{''}$
 for all $P\in {\rm Spec}\,R$ and $N(P)$ is absolutely indecomposable.
Thus $A$ is associative and separable, but not commutative, and $N$ absolutely indecomposable.
\item{} If $m=6$ then $A(P)$ is isomorphic to the direct sum of a quadratic \'{e}tale algebra and some quaternion algebra
 for all $P\in {\rm Spec}\,R$.
Thus $A$ is associative and separable, not commutative, and $N$ absolutely indecomposable.
\item{} If $m=8$ then $A(P)$ is either an octonion algebra over $k(P)$, a quaternion algebra over some quadratic field extension
of $k(P)$, or the direct sum of two quaternion algebras over $k(P)$.
 Hence $A$ is not commutative and $N$ absolutely indecomposable.
If, in particular, $A$ is not associative, then $A(P_0)$ is an octonion algebra over $k(P_0)$ for some $P_0
\in {\rm Spec}\,R$ and $N(P_0)$ is the square of its norm.
\item{} If $m=9$, then
$A(P)\cong k(P)\oplus \text{ ``some octonion algebra over } k(P)\text{''}$
for all $P\in {\rm Spec}\,R$. Thus $A$ is not associative, not commutative, and $N$ absolutely indecomposable.
\item{} If $m=10$, then
$A(P)$ is isomorphic to the direct sum of a quadratic \'{e}tale algebra and some octonion algebra over $k(P)$, or
$A(P)\cong k(P)\oplus \text{ ``some Azumaya algebra of} $ $ \text{degree 3 over } k(P)\text{''}.$
Thus $A$ is not associative, not commutative and $N$ absolutely indecomposable.
\item{} If $m=12$, then
$A(P)$ is isomorphic to the direct sum of some octonion algebra and some quaternion algebra for all $P\in {\rm Spec}\,R$.
 Again $A$ is not commutative, not associative and $N$ absolutely indecomposable.
\item{} If $m=16$ and $A$ is associative, then  $A(P)$ is an central simple algebra over $k(P)$ for all
$P\in {\rm Spec}\,R$. Hence $A$ is an Azumaya algebra over $R$ with norm $N$. If $A$ is not associative, then
$A(P_0)$ is an octonion algebra over some quadratic field extension of $k(P_0)$, for some $P_0\in {\rm Spec}\,R$.
\end{enumerate}
\end{proof}

\begin{lemma}  Let $N:A\to R$ be a nondegenerate quartic form $N$ on $A$ permitting composition.
\begin{enumerate}
\item{} Suppose $A$ can be written as the direct sum of two non-zero ideals.
 Then $A=A_1\oplus A_2$ and $N(x_1+x_2)=n_1(x_1)n_2(x_2)$. Either $A_1$ and $A_2$ are composition algebras
 and  $n_1$ and $n_2$ are their norms, or $A_1=R$, $n_1(x_1)=x_1$ and $n_2$ is a nondegenerate  cubic form permitting composition.

\item{} Suppose $A$ can be written as the direct sum of three non-zero ideals.
 Then $A\cong R\oplus R\oplus A_3$ and $N(x_1+x_2+x_3)=x_1 x_2 n_3(x_3)$ with
  $A_3$ a composition algebra over $R$ with norm $n_3$.

\item{} Suppose $A$ can be written as the direct sum of four non-zero ideals.
 Then $A\cong R\oplus R\oplus R\oplus R$ and $N(x_1+x_2+x_3+x_4)=x_1x_2x_3x_4$.
\end{enumerate}
Each $N$ is  absolutely indecomposable.
\end{lemma}

\begin{proof} \begin{enumerate}
\item{} Suppose $A$ can be written as the direct sum of two non-zero ideals.
 Then $A=A_1\oplus A_2$ and $N(x_1+x_2)=n_1(x_1)n_2(x_2)$ with $n_1$ and $n_2$ being
either two nondegenerate quadratic forms permitting composition, or $A_1=R$, $n_1(x_1)=x_1$ and
$n_2$ is a nondegenerate  cubic form permitting composition (Theorem 4 (i))
 Since $R$ is a domain, $A_1$ and $A_2$ have full support in both cases.

\item{} Suppose $A$ can be written as the direct sum of three non-zero ideals.
 Then $A\cong A_1\oplus A_2\oplus A_3$ with $A_1=A_2=R$, and $N(x_1+x_2+x_3)=n_1(x_1)n_2(x_2)n_3(x_3)$,
  $n_1(x_1)=x_1$, $n_2(x_2)=x_2$ and $n_3$ a nondegenerate quadratic form permitting composition (Theorem 4 (i)).
Thus $A_3$ is a composition algebra over $R$, since $R$ is a domain.
\end{enumerate}
 The rest of the assertion is clear.
\end{proof}

 \section{Examples}

Let $k$ be an infinite field. Let $A$ be a unital nonassociative algebra of constant rank over $R$
 which is finitely generated projective and faithful as $R$-module. Let  ${\rm Cay}(D, \mu )$
 be the classical Cayley-Dickson doubling of the composition algebra $D$ with scalar $\mu\in R^\times$ (cf.
 for instance [P]).

\subsection{} Let $R=k[t]$ be the polynomial ring over $k$.
 Suppose that $k$  has characteristic $0$ or greater than $3$ and that there exists a nondegenerate cubic form
on $A$ permitting composition. If $A$ is the direct sum of two non-zero ideals, then
 $A=R\oplus C$ and $N(x_1+x_2)=x_1n_C(x_2)$, where $C=C_0\otimes R$ is a composition algebra defined over $k$
 with norm $n_C$ [P, 6.8]. Hence both  $A=R\oplus C$ and $N$ are defined over $k$ and $N$ is absolutely
indecomposable.

If $n_{l/k}$ is the norm of a cubic field extension of $k$, then $N=n_{l/k}\otimes_k R$ is another example of an
indecomposable cubic form over $R$ permitting composition, again $N$ is defined over $k$.

\smallskip
Suppose now that $k$  has characteristic $0$ or greater than $4$ and
that there exists a nondegenerate quartic form on $A$ permitting composition.
 The following are examples of such an $(A,N)$:\\
 $(i)$ $A=R$ and $N(x)=x^4$;\\
 $(ii)$ A composition algebra $A$ over $R$ with norm $n$ and $N(x)=n(x)^2$.
 Then both $A$ and $N$ are defined over $k$ and $N$ is absolutely indecomposable.\\
$(iii)$  A separable quartic ring extension $A$ of $R$  and $N$ its norm, e.g.
$N=n_{l/k}\otimes_k R$ where $n_{l/k}$ is the norm of a separable quartic field extension of $k$.\\
$(iii)$  A composition algebra $A$ of constant rank $4$ or $8$ over its center, which is a separable quadratic ring
 extension $R'$ of $R$. $A$ has a (unique) norm $n_{A/R'}$ and $N(x)=n_{R'/R}(n_{A/R'}(x))$.
Thus $R'$  is
isomorphic to $R\oplus R$ or to $k(\sqrt{c}) \otimes_k R$.
 Say $R'=k(\sqrt{c}) \otimes_k R\cong k(\sqrt{c})[t]$.
Then either $A$ is isomorphic to ${\rm Mat}_2 (R')$ or ${\rm Zor}\,R'$, or it is without zero divisors
and defined over $k(\sqrt{c})$.
\\
$(iv)$ An Azumaya algebra $A$ of rank 16 over $R$ and $N$ its norm.\\
$(v)$ $A=C\oplus D$ and $N(x_1+x_2)=n_{C}(x_1)n_D(x_2)$ with $C$, $D$ two composition algebras over $R$ which are defined over $k$.
 $A$ and $N$ are defined over $k$ with norms $n_C$, $n_D$.
$N$ is absolutely indecomposable.\\
$(vi)$ $A=R\oplus A_2$, where $A_2$ is an algebra over $R$ with a cubic form $n_2$
permitting composition; that means, $N$ is absolutely indecomposable. For instance,
\begin{enumerate}
\item{} $A=R\oplus R$ and $N(x_1+x_2)=x_1x_2^3$;
\item{} $A=R\oplus R\oplus A_2$ and $N(x_1+x_2+x_3)=x_1x_2n(x_3)$, where $A_2$ is a cubic \'{e}tale or an
Azumaya algebra of rank 9 over $R$ and $n$ its norm;
\item{} $A=R\oplus R \oplus C$ and $N(x_1+x_2+x_3)=x_1x_2n_C(x_3)$, where $C$ is a composition algebra over $R$
 defined over $k$ with norm $n_C$. $A$ and $N$ are defined over $k$.
\end{enumerate}

$(vii)$ $A=R\oplus R\oplus C$ and $N(x_1+x_2+x_3)=x_1x_x n_C(x_3)$,
where $C$ is a composition algebra over $R$ Both $A$ and $N$ are defined over $k$, $N$ is absolutely indecomposable.

\subsection{} Let $R = k [t, \frac{1}{t}]$ be the ring of Laurent polynomials
over $k$.
  Suppose that  $k$  has characteristic $0$ or greater than $3$ and that there exists a nondegenerate cubic form on $A$ permitting composition.
   If $A$ is the direct sum of two non-zero ideals, then
$A=R\oplus C$ where $C$ is a composition algebra over $R$ with norm $n_C$ and
$N(x_1+x_2)=x_1n_C(x_2)$. Thus $N$ is absolutely indecomposable and either $C$ (hence $A$ and $N$) is defined over $k$ or
it is isomorphic to ${\rm Cay}(D, \mu t)$, where $D$ is a composition algebra without zero divisors
of half the rank of $C$ which is defined over $k$, and $\mu\in k^\times$ arbitrary [Pu1].

Other examples of absolutely indecomposable cubic forms over $R$ permitting composition are of the type $N=n_{l/k}\otimes_k R$ where
$n_{l/k}$ is the norm of a cubic field extension of $k$, or of the kind $A=J(R,\mu t)$ and
$N(u,v,w)=u^3+\mu tv^3+\mu^{-1}t^{-1} w^3-3\mu t uvw$ where $J(R,\mu t)$ is the first Tits construction starting with $R$,
$\mu\in k^\times$; i.e., $J(R,\mu t)=R(\sqrt[3]{\mu t})$, cf. Example 1.

\smallskip
Suppose now that $k$  has characteristic $0$ or greater than $4$ and that there exists a nondegenerate quartic form on $A$
 permitting composition. The following are examples of such $(A,N)$:\\
 $(i)$ $A=R$ and $N(x)=x^4$.\\
 $(ii)$  $A$ a composition algebra over $R$ with norm $n$  and $N(x)=n(x)^2$.
Then $A$ is split and isomorphic to $R\oplus R$,  ${\rm Mat}_2 (R)$ or ${\rm Zor}\,R$, or $A=A_0\otimes_k R$
with $A_0$ a composition division algebra over $k$.
 Consider the non-split case: If it has rank 2 it is either isomorphic to $k(\sqrt{c}) \otimes_k R$ for some $c\in k^\times$ which is not a square,
   or to ${\rm Cay }(R, \mu t)$ with $\mu \in
k^{\times}$. Every composition algebra of rank greater than 2 without zero divisors is either defined over $k$
or it is isomorphic to ${\rm Cay}(T, \mu t)$ where $T$ is a composition algebra without zero divisors
 of half the rank which is defined over $k$, and $\mu\in k^\times$ arbitrary. \\
$(iii)$ A separable quartic ring extension of $R$  with norm $N$.\\
$(iv)$ A composition algebra of constant rank $4$ or $8$ with center $R'$, which is a separable quadratic ring extension of $R$
with (unique) norm $n_{A/R'}$ and $N(x)=n_{R'/R}(n_{A/R'}(x))$.

 In particular, it is possible that $R'$  is
isomorphic to $k(\sqrt{c}) \otimes_k R$ or to ${\rm Cay }(R, \mu t)$ with $\mu \in
k^{\times}$ (unless $R'$ is not \'{e}tale, in that case there might be others).
 Suppose that $R'=k(\sqrt{c}) \otimes_k R\cong
k(\sqrt{c})[t,\frac{1}{t}]$ for some quadratic field extension $k(\sqrt{c})$ of $k$. Then $A$ is
isomorphic to ${\rm Mat}_2 (R')$ or ${\rm Zor}\,R'$, or it is a composition division algebra over $R'$.
 It is either defined over $k(\sqrt{c})$
or it is isomorphic to ${\rm Cay}(T, \mu t)$ where $T$ is a composition algebra without zero divisors
 of half the rank which is defined over $k(\sqrt{c})$, and $\mu\in k(\sqrt{c})^\times$ arbitrary.
  \\
$(iv)$ An Azumaya algebra of rank 16 over $R$ with norm $N$.\\
$(v)$ $A=C\oplus D$ and $N(x_1+x_2)=n_{C}(x_1)n_D(x_2)$ with $C$, $D$ two composition algebras over $R$ which
are defined over $k$.
$N$ is absolutely indecomposable.\\
$(vi)$ $A=R\oplus A_2$, where $A_2$ is an algebra over $R$ with a cubic form $n_2$
permitting composition; so $N$ is absolutely indecomposable.

$(vii)$ $A=R\oplus R\oplus C$ and $N(x_1+x_2+x_3)=x_1x_x n_C(x_3)$,
where $C$ is a composition algebra over $R$ with norm $n_C$; i.e., $C$ is defined over $k$
or isomorphic to ${\rm Cay}(T, \mu t)$ where $T$ is a composition algebra without zero divisors
 of half the rank which is defined over $k$,  $\mu\in k^\times$ arbitrary. $N$ is absolutely indecomposable.

\section{Forms permitting composition over locally ringed spaces}

 Let $X$ be a locally ringed space with structure sheaf $\mathcal{O}_X$ such that $d!\in H^0(X,\mathcal{O}_X^\times)$.

 An $\mathcal{O}_X$-algebra $\mathcal{A}$ is called {\it alternative} if $x^2y=x(xy)$ and
$yx^2=(yx)x$ for all sections $x,y$ of $ \mathcal{A}$ over the same open subset of $X$.
 An associative  $\mathcal{O}_X$-algebra $\mathcal{A}$ is called an {\it Azumaya algebra} if
$\mathcal{A}_P\otimes_{\mathcal{O}_{P,X}} k(P)$ is a central simple algebra over $k(P)$ for all $P\in X$.

Let $ \mathcal{A}$ be an $\mathcal{O}_X$-algebra together with a nondegenerate form $N:\mathcal{A}\to \mathcal{O}_X$
of degree $d$ permitting composition; i.e.,
 $N(xy)=N(x)N(y)$
for all sections $x,y$ of $\mathcal{A}$ over the same open subset of $X$.
Let $1=1_\mathcal{A}\in H^0(X, \mathcal{A})$ be the unit element of $\mathcal{A}$.
Suppose always that $N(1)=1$. Then $\mathcal{A}$ has full support.

\begin{remark} (i) Let $\mathcal{A}=\mathcal{A}_1\oplus \mathcal{A}_2$ be the direct sum of two non-zero ideals
of $ \mathcal{A}$. Since $\mathcal{A}$ is locally free of finite rank as  $\mathcal{O}_X$-module by our convention 2.4,
 so are $\mathcal{A}_1$ and $\mathcal{A}_2$.
\\(ii) Let $X$ be a scheme over the affine scheme $Y={\rm Spec}\,R$, $H^0(X,\mathcal{O}_X)=R$.  If $\mathcal{A}$
is globally free as an $\mathcal{O}_X$-module then $(\mathcal{A},N)$ is defined over $R$.
(The proof is analogous to the one of [P, 1.10].)
\end{remark}

Let $\theta : \mathcal{A} \times \dots \times \mathcal{A}\to \mathcal{O}_X$ be the $d$-linear form associated with $N$.
For $i=1,\dots,d$  define a form $T_i:\mathcal{A}\to \mathcal{O}_X$ of degree $i$ via
$$T_i(x)=\binom{n}{i}\theta(x,\dots, x,1,\dots,1)\,\,\text{  ($i$-times }x).$$ Then
$$N(x)=T_d(x) \text{ and } T_1(x)=d\theta(x,1,\dots,1)$$
for all sections $x$ of $\mathcal{A}$ over the same open subset of $X$.

The form $T:\mathcal{A}\to \mathcal{O}_X$, $T(x)= T_1(x)$ is called the {\it trace}.
Put $T_0(x)=1$ and $T_{d+q}(x)=0$ for $q>0$, then $T(a1)=da\theta(1,\dots,1)$ for all $a$ in $\mathcal{O}_X$.
Define $\mathcal{A}_0={\rm ker}\, T$. $\mathcal{A}$ is the direct sum of $\mathcal{O}_X$-modules
$\mathcal{A}=\mathcal{O}_X1_\A\oplus \mathcal{A}_0.$

We assume from now on that $d\geq 3$.
 Our results from Section 3 easily adapt to the setting of locally ringed spaces:  $\mathcal{A}$ is alternative and
 $B: \mathcal{A}\times  \mathcal{A}\to  \mathcal{O}_X$, with $B(x,y)=T(xy)$ for all sections $x,y$ of $ \mathcal{A}$ over
 the same open subset of $X$, is a nondegenerate symmetric  bilinear form on $ \mathcal{A}$ which is associative.
Every section $x$ of $\mathcal{A}$ over the same open subset of $X$ satisfies
$$x^d-T_1(x)x^{d-1}+T_2(x)x^{d-2}-\cdots+(-1)^d T_d(x)1=0.$$

 From now on let $X$ be an integral scheme.
If $\mathcal{A}=\mathcal{A}_1\oplus \cdots \oplus \mathcal{A}_r$ with $\mathcal{A}_i$
non-zero ideals of $\mathcal{A}$, then
 $$N(x)=N_1(x_1)\cdots N_2(x_r)$$
where $x=x_1+\dots +x_r$, $x_i\in \mathcal{A}_i$ for all $i$, and each $N_i$ is a nondegenerate
form of degree $d_i$ on $\mathcal{A}_i$, $d=d_1+\dots+d_r$, which permits composition.

\begin{lemma} (i) $\mathcal{A}(P)$ is separable for all $P\in X$.
\\(ii) If $\mathcal{A}$ is not associative, but $\mathcal{A}(P)$ is simple for all $P\in X$,
  then there is at least one $P\in X$ such that
 $\mathcal{A}(P)$ is an octonion algebra over some separable field extension of $k(P)$.
\end{lemma}

\begin{proof} (i) We know that for all $P\in X$, $ \mathcal{A}(P)$ is a non-zero $k(P)$-algebra
 such that there exists a nondegenerate form $N(P)$ of degree $d$ on $\mathcal{A}(P)$ permitting composition
 (2.3 and Lemma 1).
Thus $ \mathcal{A}(P)$ is a finite dimensional separable $k(P)$-algebra for all $P\in X$.
\\ (ii) follows immediately.
\end{proof}

\begin{proposition} (i) If $N$ is a cubic form then $\mathcal{A}$ has rank $1$, $2$, $3$, $5$ or $9$.\\
(ii) If $N$ is a quartic form, then $\mathcal{A}$ has rank $1$, $2$, $3$, $4$, $5$, $6$, $8$, $9$, $10$, $12$ or $16$.
\end{proposition}

\begin{proof} For all $P\in {\rm Spec}R$ the residue class algebra $\mathcal{A}(p)=\mathcal{A}_P\otimes k(P)$ is  a nonzero alternative algebra
 over the residue class field $k(P)$ together with a nondegenerate form $N(P):\mathcal{A}(P)\to k(P)$
 which permits composition  (2.3 and Lemma 1).
 This implies that if $N$ is a cubic form then $\mathcal{A}(P)$ must have rank
$1$, $2$, $3$, $5$ or $9$, and if it is a quartic form then $\mathcal{A}(P)$ must have rank
$1$, $2$, $3$, $4$, $5$, $6$, $8$, $9$, $10$, $12$ or $16$ [S1, p.~140].
\end{proof}

We now turn to cubic forms  permitting composition.

\begin{lemma} A  nondegenerate cubic form $N$ permitting composition
 cannot be defined\\
(i) on a composition algebra over $X$  of constant rank greater than $2$;\\
(ii) on a non-split quadratic \'etale algebra $\mathcal{A}$ over $X$, unless $\mathcal{A}(P)$ is split for all $P\in X$.
\end{lemma}

 This can be proved by assuming such a form exists and using Lemma 3 locally to obtain a contradiction.

\begin{proposition} Suppose that there exists a nondegenerate cubic form $N$ on $\mathcal{A}$ permitting composition.
Then $(\mathcal{A},N)$ satisfies one of the following:\\
 (i) $\mathcal{A}=\mathcal{O}_X$ and $N(x)=x^3$;\\
 (ii) $\mathcal{A}$ is  commutative and associative of rank 2 or 3, and $N$ is absolutely
 indecomposable; if $\mathcal{A}$ is a non-split quadratic \'etale algebra over $X$, then $\mathcal{A}(P)$ must be split for all $P\in X$;
 if  $\mathcal{A}$ has rank 3 then  $\mathcal{A}(P)$ is a cubic  \'etale algebra over $k(P)$ for all $P\in X$.\\
(iii) $\mathcal{A}$ is an associative and not commutative $\mathcal{O}_X$-algebra of rank $5$, and $N$ is absolutely indecomposable.\\
(iv)   $\mathcal{A}$ is an Azumaya algebra over $X$ of  rank $9$, and $N$ is its
reduced norm.\\
(v) $\mathcal{A}$ has rank $9$ and is neither associative nor commutative, $N$ is absolutely indecomposable.
\end{proposition}

This follows from Proposition 2.

\begin{example} There is the following first Tits construction starting with the structure sheaf of
an integral scheme $X$
[Ach, 2.25]:
Let $\lb\in \, _3 {\rm Pic} X$ and $N:\lb\to \sts$ a nondegenerate cubic form, let $\check{\lb}=
{\rm Hom}_{\sts}(\lb,\sts)$ be its dual and $\langle w,\check{w}\rangle=\check{w}(w)$ the canonical pairing $\lb
\times\check{\lb}\to \sts.$ There exists a uniquely determined cubic norm $\check{N}:\check{\lb}\to \sts$ and
uniquely determined adjoints $\sharp:\lb\to \check{\lb}$ and $\check{\sharp}:\check{\lb}\to \lb$ such that
the identities listed in Example 1 (with $N$ instead of $N_L$) hold
\ignore{\begin{enumerate}
\item $\langle w,w^\sharp\rangle=N(w)1$;
\item $\langle \check{w}^{\check{\sharp}},\check{w}\rangle=\check{N}(\check{w})1$;
\item $w^{\sharp \, \check{\sharp}}=N(w)w$;
\end{enumerate}
for $w$ in $\lb$, $\check{w}$ in $\check{\lb}$. Moreover,
\begin{enumerate}
\item $\check{w}^{\check{\sharp}\,\sharp }=\check{N}(\check{w})\check{w}$;
\item $\langle w, \check{w}\rangle^2= \langle \check{w}^{\check{\sharp}}, w^{\sharp}\rangle$;
\item $\langle w,\check{w}\rangle^3=N(w)\check{N}(\check{w})$;
\item $\langle w,\check{w}\rangle w= 3\langle w,\check{w}\rangle w-w^\sharp\check{\times}\check{w}$
\end{enumerate}}
for $w,w'$ in $\lb$, $\check{w}$ in $\check{\lb}$ [Ach, 2.13].
 Define
$$\begin{array}{l}
\widetilde{\J}=\mathcal{O}_X\oplus \lb\oplus\check{\lb},\\
\widetilde{1}=(1,0,0)\in H^0(X,\J),\\
\widetilde{N}(a,w,\check{w})=a^3+N(w)+\check{N}(\check{w})-3a\langle w,\check{w}\rangle\\
(a,w,\check{w})^{\widetilde{\sharp}}=(a^2-\langle w,\check{w}\rangle, \check{w}^{\check{\sharp}}-
aw, w^\sharp-\check{w}a)
\end{array}$$
for $a\in\sts$, $w\in\lb$, $\check{w}\in\check{\lb}$, then $(\widetilde{N},\widetilde{\sharp},\widetilde{1})$
is a cubic form with adjoint and base point on $\widetilde{\J}$ and has trace form
$$\widetilde{T}((a,w,\check{w}),(c,v,\check{v}))=3ac+3\langle w, \check{v}\rangle+3\langle v,
\check{w}\rangle.$$
The induced (commutative associative) Jordan algebra $J(\widetilde{N},\widetilde{\sharp},\widetilde{1})$ is denoted by $J(\A,\Pe,N)$.
This construction yields examples of commutative associative algebras
 $\A=
 \J (\sts, \lb,  N)$ of constant rank 3 which admit a nondegenerate cubic form $\widetilde{N}$ permitting composition.
 $\A(P)$ is a cubic \'{e}tale algebra over $k(P)$ for all $P\in {\rm Spec}\,X$.
\end{example}

\begin{example} [Pu2] Let $k$ be a field of characteristic 0. Let $X$ be an elliptic curve $X$ over $k$.
 Let $\N_i$ denote a line bundle of order 3 on $X$ with $\N_0= \mathcal{O}_{X} $.  We have
 $_3 \Pic (X) = \{\N_i\,|\, 0\leq 1\leq m\}$ for some even integer $0\leq m\leq 8$.
 Every first Tits construction  over $X$ starting with
 $\mathcal{O}_X$ is isomorphic to $\A=\J (\mathcal{O}_X, \mathcal{N}_i, N_i)$ where
$N_i$ is a nondegenerate cubic form on $\N_i$.
By the Theorem of Krull-Schmidt, if $\N_i\not\cong\N_j$ and $\N_i\not\cong\N_j^\vee$
 then $\J(\sts,\N_i, N_i)\not \cong\J(\sts,\N_j, N_j)$. $\J(\sts,\N_0, N_0)$ is defined over $k$.
\end{example}

\section{Curves of genus zero}

\begin{lemma} Let $X$ be a curve of genus zero over a field $k$ of characteristic not 2 or 3.
Let $\mathcal{A}$ be an algebra over $X$ of rank 2 carrying a nondegenerate form $N$ permitting composition.
If $N$ is a cubic or quartic form then $(\mathcal{A},N)$ is defined over $k$.
\end{lemma}

\begin{proof} $\mathcal{A}$, together with the nondegenerate symmetric
bilinear form $T_\mathcal{A}$, is a nondegenerate bilinear space over $X$.
If $X$ is rational, [P, 5.4] shows that $\mathcal{A}$ decomposes into the orthogonal sum of $\mathcal{O}_X$-modules
of the kind $\mathcal{O}_X(m_i)\oplus\mathcal{O}_X(-m_i)$ for suitable $m_i>0$.
 Hence, if $\mathcal{A}$ has rank two, it must be globally free as an $\mathcal{O}_X$-module.

If $X$ is nonrational, there is a field extension $k'/k$ such that $X'=X\times_kk'$
becomes rational. If $\mathcal{A}$ is an algebra as in our assumption then so is $\mathcal{A}\otimes \mathcal{O}_{X'}$
and since this is globally free, so is $\mathcal{A}$.  By Remark 7 (ii), in both cases $(A,N)$ is defined over $k$.
\end{proof}

\begin{lemma}  Let $X$ be a curve of genus zero over a field $k$ of characteristic not 2 or 3.
 Every first Tits construction over $X$ starting with $\mathcal{O}_X$
 is defined over $k$.
\end{lemma}

\begin{proof} Every first Tits construction  over $X$ starting with $\mathcal{O}_X$ is of the kind
$J (\mathcal{O}_X, \mathcal{L}, N_{\mathcal{L}})$, where
$\mathcal{L} \in \, _3 {\rm Pic} X$ (Example 2). However,
 ${\rm Pic}\,X\cong\mathbb{Z}$, so we only have $J (\mathcal{O}_X,\mathcal{O}_X, \mu)\cong
 J(k,\mu)\otimes_k \sts$, $\mu\in k^\times$.
\end{proof}

\subsection{} Let $k$ be a field of characteristic $0$ and $X$ be a curve of genus zero over $k$.
 If $X$ is not rational, let $D_0=(a,b)_k$ be the quaternion division algebra associated with $X$.
 Let $\mathcal{A}$ be an algebra (automatically of constant rank) over $X$ such
that there exists a nondegenerate cubic form $N$ on $\mathcal{A}$ permitting composition.
\begin{enumerate}
 \item If $\mathcal{A}$ has rank 1, then $\mathcal{A}=\mathcal{O}_X$ and $N(x)=x^3$.
\item If $\mathcal{A}$ has rank 3 it  is commutative associative and $N$ is absolutely indecomposable.
For instance, let $k'$ be a cubic field extension of $k$, then $\mathcal{A}=k'\otimes_k \mathcal{O}_X$ carries a nondegenerate
cubic form permitting composition.
  \item If $\mathcal{A}$ is the direct sum of two non-zero ideals $\mathcal{A}=\mathcal{A}_1\oplus \mathcal{A}_2$, then
$N(x_1+x_2)=N_1(x_1)N_2(x_2)$ and  $N_1$ must be a nondegenerate linear form  and $N_2$ a nondegenerate
 quadratic one (or vice versa). It follows that
 $\mathcal{A}_1=\mathcal{O}_X$, $N_1=id$, and $\mathcal{A}_2$ is an algebra of degree 2 over $X$ of
 constant rank with a nondegenerate quadratic form $N_2$ permitting composition.
Therefore $\mathcal{A}_2$ is a composition algebra over $X$ of constant rank with norm $N_2$.
 Hence
 $$\mathcal{A}=\mathcal{O}_X\oplus \mathcal{C}$$
 with $C$ a quadratic \'etale algebra, a quaternion or an octonion algebra over $X$, and $N$ is absolutely indecomposable.
  (Note that $\mathcal{A}$ cannot be the direct sum of more than two non-zero ideals.)
   By  [P, 4.4], one of the following holds:\\
 (i) $\mathcal{C}$ (and thus $\mathcal{A}$) is defined over $k$.\\
 (ii) $\mathcal{C}$ is a split quaternion or octonion algebra.\\
 (iii) $X$ is not rational and $\mathcal{C}\cong {\rm Cay} (\mathcal{D},\mathcal{P}, N_\mathcal{P})$,
 where $\mathcal{D}=D_0\otimes \mathcal{O}_X$, $\mathcal{P}$ is a locally free right
$\mathcal{D}$-module of rank one and norm one, and $N_\mathcal{P}$ is a norm on it.
More precisely, let $\E_0$ be the indecomposable $\sts$-module of rank 2
 described in [P, 4.3]. Then we know that $\mathcal{P}=\mathcal{P}_1\oplus \mathcal{P}_2$ with
$\mathcal{P}_1=\mathcal{L}(mP_0)\otimes \check{{\mathcal{E}}}_0$ and
$\mathcal{P}_2=\mathcal{L}((-m+1)P_0)\otimes \check{\mathcal{E}}_0$ for some integer $m\geq 0$ uniquely
determined by $\mathcal{C}$, where $P_0$ is a closed point of $X$ of minimal degree.
  \item If $\mathcal{A}$ is the direct sum of three non-zero ideals
  $\mathcal{A}=\mathcal{A}_1\oplus \mathcal{A}_2\oplus \mathcal{A}_3$, then
$N(x_1+x_2+x_3)=N_1(x_1)N_2(x_2)N_2(x_3)$ and  $N_i$ must be a nondegenerate linear form ($i=1,2,3$)
 It follows that $N(x_1+x_2+x_3)=x_1x_2x_3$ and $\mathcal{A}$ is defined over $k$.

 \item Suppose that $\mathcal{A}$ does not decompose into the direct sum of non-zero ideals.
\\ If $\mathcal{A}$ has rank 9 and  is associative then  $\mathcal{A}$ is an Azumaya algebra over $X$ of
rank $9$ and $N$ its reduced norm. If $X$ is rational, then either we have
$$\mathcal{A}\cong \mathcal{E}nd_X(\mathcal{O}_X(m_1)\oplus \mathcal{O}_X(m_2)\oplus\mathcal{O}_X(m_3))$$
with $m_i\in\mathbb{Z}$, or both $\mathcal{A}$ and $N$ are defined over $k$ and we have
$$\mathcal{A}\cong
\sigma^*D$$
 with $D$ a central simple division
 algebra over $k$ ([Kn1, VII (3.1.1)], [Kn2], see also [Ach]).
If $X$ is nonrational, we can give some examples
 of Azumaya algebras of rank 9 over $X$: \\
 For instance, again the trivial case that $\mathcal{A}\cong
\sigma^*D$
 with $D$ a central simple division algebra over $k$.\\
  For every locally free $\sts$-module $\mathcal{E}$
 of constant rank 3, $\E nd_X(\E)$ is an Azumaya algebra of rank 9 where
 we have the following possibilities for $\E$:
$$\E=\lb(m_1P_0) \otimes\E_0\oplus\lb(m_2P_0) \text{ and }
\E=\lb(n_1P_0)\oplus\lb(n_2P_0)\oplus\lb(n_3P_0)$$
 with $m_1,m_2,n_1,n_2,n_3\in\mathbb{Z}$. Hence
\\
(i)  \[\E nd_X(\E)=
\left [\begin {array}{ccc}
\E nd(\E_0)&\lb(-aP_0)\otimes\E_0\\
\noalign{\smallskip}
\lb(aP_0)\otimes\check{\E_0}&\sts\\
\end {array}\right ]
\]
(ii) \[ \E nd_X(\E)=
\left [\begin {array}{ccc}
\sts & \lb(cP_0) & \lb(bP_0)\\
\noalign{\smallskip}
\lb(-cP_0) &  \sts & \lb((b-c)P_0)  \\
\noalign{\smallskip}
\lb(-bP_0) & \lb((c-b)P_0) &  \sts \\
\end {array}\right ]
\]
with $a=m_2-m_1$, $b=n_1-n_3$ and $c=n_1-n_2$  [Ach, 4.4].
\end{enumerate}

\section{Degenerate forms permitting composition}

 \subsection{} Let us consider forms $N:A\to R$ permitting composition over rings $R$ with $d!\in R^\times$,
where $A$ is a unital nonassociative $R$-algebra which is finitely generated projective
as $R$-module. We now look at the case where $N:A\to R$ is degenerate; i.e., where the {\it radical}
$${\rm rad}\,N=\{x\in A\,|\,\theta(x,x_2,\dots,x_d)=0 \text{ for all }x_i\in A\}$$
 is non-zero. As before, let $B:A\times A \to R$, $B(x,y)=T(xy)$ where $T$ is the trace of $A$.

\begin{lemma}  Let $A$ be an $R$-algebra together with a form $N$ of degree $d$ on
$A$ permitting composition. Let $D$ be a subalgebra of $A$ which is  maximal  among all subalgebras $E$ of $A$
which have a nondegenerate restriction $N|_E$. Then $(A,B)=(D,B|_D)\perp (D^\perp, B|_{D^\perp})$
and $DD^\perp\subset D^\perp$ as well as $D^\perp D\subset D^\perp$.\\
\end{lemma}

\begin{proof}
Since $N_D$ is nondegenerate, so is $B|_D$. Thus
$(A,B)=(D,B|_D)\perp (D^\perp, B|_{D^\perp})$ by [Kn1, I(3.6.2)]. Let $x\in D$ and $y\in D^\perp$.
Then $B(z,xy)=B(zx,y)=0$ and $B(yx,z)=B(y,xz)=0$ for all $z\in D$, thus $xy\in D^\perp$ and $yx\in D^\perp$.
Since $B$ is an associative symmetric bilinear form, ${\rm rad}\,B$ is a two-sided ideal. We have
${\rm rad}(B)={\rm rad}(N)$  by Lemma 1 (i).
\end{proof}

For $d=2$ and $R$ a field this was proved in [K-S, 1.2].

\begin{remark}  The radical ${\rm rad}\,N$ of a cubic form $N$ permitting composition on an $R$-algebra $A$ is a two-sided
ideal [B-B, Lemma 2] (indeed, this is true for the radical of any  form of  degree $d$ permitting composition [B]).
If $A$ is also an algebra of degree 3 as defined in [B-B], then ${\rm rad}\,N$ is a nilideal. If, additionally, $R$ does not contain any non-zero nilpotent elements,
then ${\rm rad}\,N$ is the maximal nilideal of $A$, that means the radical [B-B, Lemma 5].
\end{remark}

 For quadratic forms permitting composition, the radical can be annihilated by a suitable exponent
which depends on the dimension of $A$ ([K-S] or [M2]). For degenerate cubic forms permitting composition of the kind
$N(a,x)=aN_C(x)$ where $N_C$ is a quadratic form permitting composition,
 the radical can be annihilated by exactly that exponent
which depends on the dimension of $C$, since in that case ${\rm rad}\,N_0=0\oplus {\rm rad}\,N_C$.

\subsection{} Let $X=\mathbb{P}_R^n$ be the $n$-dimensional projective space over $R$, that is $X={\rm Proj}\,S$
where $S=R[t_0,\dots,t_n]$ is the polynomial ring in $n+1$ variables over $R$, equipped
 with the canonical grading $S=\oplus_{m\geq 0}S_m$.  We have ${\rm rank}\, S_m=\binom{m+n}{n}$.
We know that $\mathcal{O}_X(m)$ is a locally free $\mathcal{O}_X$-module of rank one
for each $m\in \mathbb{Z}$ and
$$H^0(X,\mathcal{O}_X(m))=S_m \text{ for }m\geq 0,$$
$$H^0(X,\mathcal{O}_X(m))=0 \text{ for }m< 0.$$

\begin{example}  Let $\mathcal{C}$ be the split octonion algebra
$${\rm Zor}(\mathcal{O}_X(l)\oplus \mathcal{O}_X(m)
\oplus \mathcal{O}_X(-l-m),\alpha)$$
 over $X$  with norm $n_C$ as defined in [P, 3.3] ($l,m$ positive integers).
Let $\mathcal{A}=\mathcal{O}_X\oplus \mathcal{C}$ and $N((x_1,x_2))=x_1n_C(x_2)$ for all sections $x_1$ in $\mathcal{O}_X$,
$x_2$ in $\mathcal{C}$. Then $N(1)=1$, $N$ permits composition, and $N$ is absolutely indecomposable.
We get
\[ A=H^0(X, \mathcal{A})=R\oplus H^0(X, \mathcal{C})=R\oplus
\left [\begin {array}{cc}
R&S_l\oplus S_m\\
S_{l+m}&R\\
\end {array}\right ]
\]

 \smallskip\noindent
with the algebra multiplication in $H^0(X, \mathcal{C})$ as described in [P, 3.8].
$A$ is an alternative $R$-subalgebra  of $S\oplus {\rm Zor}\,S$ of rank
$$3+\binom{l+n}{n}+\binom{m+n}{n}+\binom{(l+m)+n}{n}$$
 and $x^3-T_1(x)x^2+T_2(x)x-T_3(x)1=0$
for each $x\in A$. If $n=1$ then ${\rm rank}_RH^0(X, \mathcal{A})=
6+2(l+m)\geq 10$
must be even. If $R$ is a field then the cubic form $N_0=N(X)$ restricted to the subalgebra
\[
R\oplus
\left [\begin {array}{cc}
R&0\\
\noalign{\smallskip}
0&R\\
\end {array}\right ]
\]

 \smallskip\noindent
 of rank $3$ is nondegenerate and

\[ {\rm rad}\,N_0=0\oplus {\rm rad}\,(H^0(X, \mathcal{C}))=
0\oplus
\left [\begin {array}{cc}
0&S_l\oplus S_m\\
\noalign{\smallskip}
S_{l+m}&0\\
\end {array}\right ]
\]

 \smallskip\noindent
 is the radical of $A$ [P, 3.8]. We have

 \[ ({\rm rad}\,N_0)^2=
0\oplus
\left [\begin {array}{cc}
0&0\\
\noalign{\smallskip}
S_{l+m}&0\\
\end {array}\right ]
\]

and $({\rm rad}\,N_0)^3=0$.
\end{example}

\begin{example} Let $\mathcal{F}=\mathcal{O}_X(m_1)\oplus\mathcal{O}_X(m_2)\oplus\mathcal{O}_X(m_3)$,
then  $\mathcal{A}=\mathcal{E}nd_X(\mathcal{F})$ is an Azumaya algebra over $X$ of constant rank $9$. We have

\smallskip
\[\mathcal{A}=
\left [\begin {array}{ccc}
\mathcal{O}_X&\mathcal{O}_X(a)&\mathcal{O}_X(b)\\
\noalign{\smallskip}
\mathcal{O}_X(-a)&\mathcal{O}_X&\mathcal{O}_X(b-a)\\
\noalign{\smallskip}
\mathcal{O}_X(-b)&\mathcal{O}_X(a-b)&\mathcal{O}_X
\end {array}\right ]
\]

 \smallskip\noindent
 with $a=m_1-m_2$, $b=m_1-m_3$, the right hand side being equipped with the usual matrix multiplication.
$H^0(X, \mathcal{A})$ is a unital associative $R$-algebra of degree $3$ which admits
a cubic form $N_0=N(X):H^0(X, \mathcal{A})\to H^0(X, \mathcal{O}_X)$ permitting composition which satisfies $N_0(1_A)
=1$.
 \begin{enumerate}
 \item{} If $a,b>0$ and $b-a > 0$ then

\smallskip
\[ H^0(X, \mathcal{A})=
\left [\begin {array}{ccc}
R&S_a&S_b\\
\noalign{\smallskip}
0&R&S_{b-a}\\
\noalign{\smallskip}
0&0&R
\end {array}\right ]
\]
 has rank $$3+\binom{a+n}{n}+\binom{b+n}{n}+\binom{(b-a)+n}{n}.$$
In particular, if $n=1$ then ${\rm rank}_RH^0(X, \mathcal{A})=
6+2b\geq 8$ must be even.
If $R$ is a field then

\smallskip
\[{\rm rad}\,N_0=
\left [\begin {array}{ccc}
0&S_a&S_b\\
\noalign{\smallskip}
0&0&S_{b-a}\\
\noalign{\smallskip}
0&0&0
\end {array}\right ]
\]

 is the radical of $N_0$, i.e. the radical of $A$. We get

\smallskip
\[ ({\rm rad}\,N_0)^2=
\left [\begin {array}{ccc}
0&0&S_b\\
\noalign{\smallskip}
0&0&0\\
\noalign{\smallskip}
0&0&0
\end {array}\right ]
\]
and $({\rm rad}\,N_0)^3=0$.

 \item{} If $ a=b>0$ then

\smallskip
\[ H^0(X, \mathcal{A})=
\left [\begin {array}{ccc}
R&S_a&S_a\\
\noalign{\smallskip}
0&R&R\\
\noalign{\smallskip}
0&R&R
\end {array}\right ]
\]

\smallskip\noindent
has odd rank $5+2\binom{a+n}{n}$. In particular, if $n=1$ then ${\rm rank}_R H^0(X, \mathcal{A})=5+2(a+1)\geq 9$.
If $R$ is a field then

\smallskip
\[ {\rm rad}\,N_0=
\left [\begin {array}{ccc}
0&S_a&S_a\\
\noalign{\smallskip}
0&0&0\\
\noalign{\smallskip}
0&0&0
\end {array}\right ]
\]
and $({\rm rad}\,N_0)^2=0$. Again,  ${\rm rad}\,N_0$ is the radical of $A$.

 \item{}  If $a=0$ and $b>  0$ then

\smallskip
\[ H^0(X, \mathcal{A})=
\left [\begin {array}{ccc}
R&R&S_b\\
\noalign{\smallskip}
R&R&S_b\\
\noalign{\smallskip}
0&0&R
\end {array}\right ]
\]

\smallskip\noindent
has odd rank $5+2\binom{b+n}{n}$. In particular, if $n=1$ then
${\rm rank}_R H^0(X, \mathcal{A})=5+2(b+1)\geq 9$.
If $R$ is a field then

\[ {\rm rad}\,N_0=
\left [\begin {array}{ccc}
0&0&S_b\\
\noalign{\smallskip}
0&0&S_b\\
\noalign{\smallskip}
0&0&0
\end {array}\right ]
\]
and $({\rm rad}\,N_0)^2=0$.
\end{enumerate}
\end{example}

Analogously, one can construct classes of degenerate forms of degree higher than 3 over $R$ which permit
composition.

\smallskip
{\it Acknowledgements:} The author would like to acknowledge the support of the
``Georg-Thieme-Ged\"{a}chtnisstiftung'' (Deutsche Forschungsgemeinschaft)
 during her stay at the University of Trento,
 and to thank the Department of Mathematics of the University of Trento for its hospitality.


\noindent
\begin{thebibliography}{[B-C-R]}

\bibitem[Ach] {} Achhammer, G.,  {\it Albert Algebren \"{u}ber lokal geringten R\"{a}umen}. PhD Thesis, FernUniversit\"{a}t
Hagen, 1995.

\bibitem[A] {} Albert, A. A., {\it Quadratic forms permitting composition}. Ann. of Math. (2) (1942), 161-177.

\bibitem[B] {} Bergmann, A.,  {\it Formen auf Moduln \"{u}ber kommutativen Ringen beliebiger Charakteristik}.
J. Reine Angew. Math. 219 (1965), 113-156.

\bibitem[B-B] {} Baumgartner, E., Bergmann, A., {\it Nichtausgeartete Kompositionsalgebren vom Grad 3}.
J. Reine Angew. Math.  268/269 (1974), 324-337.

\bibitem[H]{7} Hartshorne, R.,  ``Algebraic geometry''. Graduate Texts in
Mathematics, vol. 52, Springer-Verlag, Berlin-Heidelberg-New York, 1977.

\bibitem[K]{11} Keet, A., {\it Decomposition of a higher degree form}. Comm. Alg. 30 (10) (2002), 4945-4963.

\bibitem[Kn1]{11} Knus, M.-A.,  ``Quadratic and hermitian forms over rings''.
Springer-Verlag, Berlin-Heidelberg-New York,
1991.

\bibitem[Kn2]{11} Knus, M.-A., {\it Quaternionic modules over $\mathbb{P}^2(\mathbb{R})$}. Proceedings of the conference
``Brauer groups in ring theory and algebraic geometry'', Antwerpen 1981, eds. F. van Oystaeyen and A. Verschoren,
Lecture Notes in Math. 917, Springer-Verlag, New York - Heidelberg - Berlin, 1982.

\bibitem[K-S]{11} Kunze, R. A., Scheinberg, S., {\it Alternative algebras having scalar involutions}.
Pac. J. Math. 124 (1) (1986), 159-172.


\bibitem[L-L]{L} Legrand, S., Legrand, D., {\it Alg\'{e}bres sur un anneau munies d'une forme multiplicative
de degr\'e $n$}. J. Reine Angew. Math. 299/300 (1978), 171-184.

\bibitem[M1] {} McCrimmon, K., {\it A proof of Schafer's conjecture for infinite-dimensional forms permitting
composition}. J. Algebra 5 (1967), 72-83.

\bibitem[M2] {} McCrimmon, K., {\it Nonassociative algebras with scalar
involution}. Pacific J. Math. 116 (1985), 85-108.

\bibitem[P] {} Petersson, H.P., {\it Composition algebras over algebraic curves of genus
zero}. Trans. Amer. Math. Soc. 337(1) (1993), 473-491.

\bibitem[P-R]{PR} Petersson, H. P., Racine, M. L., {\it Jordan algebras of degree 3 and the Tits process}.
J. Algebra 98 (1986) (1), 211-243.

\bibitem[Pr]{B} Pr\'oszy\'nski, A., {\it On orthogonal decomposition of
homogenous polynomials}. Fundamenta Math. 98 (3) (1978), 201-217.

\bibitem[Pu1]{Pu} Pumpl\"un, S., {\it Composition algebras over the Laurent polynomials}.
Comm. Alg. 25(1) (1997), 229-233.

\bibitem[Pu2]{Pu}   Pumpl\"un, S.,  {\it Albert algebras over curves of genus zero and one.} Preprint, 2007.
\bibitem[R] {} Roby, N., {\it Lois polyn\^{o}mes et lois formelles en th\'{e}orie des modules}. Ann. Sci. Ecole Norm. Sup.
$3^e$ ser. t. 80 (1963), 213- 348.

\bibitem[S1]{S1} Schafer, R.D., {\it Forms permitting composition}. Adv.
Math. 4 (1970), 127-148.

\bibitem[S2]{S2} Schafer, R.D., {\it On cubic forms  permitting composition}. Proc. Amer. Math. Soc. 10 (1959), 917-925.

\bibitem[S3]{S3} Schafer, R.D.,
{\it On forms of degree $n$ permitting composition}. J. of Mathematics and Mechanics (1963) 12,
777-792.

\bibitem[S4]{} Schafer, R.D., ''An Introduction to Nonassociative
Algebras''. Academic Press, Inc., New York, 1966 (reprinted by Dover Publications, Inc., New York, 1995).

\bibitem[Sl1]{} Slater, M., {\it Alternative rings with d.c.c. I.} J. Algebra 11 (1969), 102-110.

\bibitem[Sl2]{} Slater, M., {\it The socle of an alternative ring}. J. Algebra 14 (1970), 443-463.

\end{thebibliography}
\end{document}